\documentclass{article}

\usepackage{amsmath, amsthm, amsfonts,amssymb}
\usepackage{color}

\usepackage{hyperref}

\usepackage{longtable}

\usepackage{csquotes}

\usepackage{mdframed}
\theoremstyle{definition}
\newmdtheoremenv{boxProb}{Problem}
\newmdtheoremenv{boxDef}{Definition}
\newmdtheoremenv{boxCor}{Corollary}
\newmdtheoremenv{boxThm}{Theorem}
\newmdtheoremenv{compjob}{Computational Job}
\newmdtheoremenv{reqi}{Requirement}

\usepackage{graphicx}
\usepackage{subfigure}          %% Unterabbildungen

\usepackage{placeins}

\usepackage{tabularx}

\usepackage{lmodern,textcomp}

\usepackage{algorithm}
\usepackage{algpseudocode}

\usepackage{enumitem}

\usepackage{xspace} 
\newcommand\largeparbreak{\par\bigskip}

\newcommand{\inv}{^{-1}\xspace}

\newcommand{\h}{^\textsf{H}\xspace}

\newcommand{\away}[1]{}

\newcommand{\N}{\mathbb{N}\xspace}
\newcommand{\C}{\mathbb{C}\xspace}

\newcommand{\cO}{\mathcal{O}\xspace}

\newcommand{\bA}{\mathbf{A}\xspace}
\newcommand{\bE}{\mathbf{E}\xspace}

\newcommand{\bB}{\mathbf{B}\xspace}
\newcommand{\bG}{\mathbf{G}\xspace}

\newcommand{\bU}{\mathbf{U}\xspace}
\newcommand{\bV}{\mathbf{V}\xspace}
\newcommand{\bZ}{\mathbf{Z}\xspace}
\newcommand{\bD}{\mathbf{D}\xspace}

\newcommand{\bQ}{\mathbf{Q}\xspace}

\newcommand{\bF}{\mathbf{F}\xspace}
\newcommand{\bx}{\mathbf{x}\xspace}

\newcommand{\by}{\mathbf{y}\xspace}

\newcommand{\bei}[1]{{\mathbf{e}}\xspace}

\newcommand{\bM}{\mathbf{M}\xspace}

\newcommand{\bC}{\mathbf{C}\xspace}

\newcommand{\bI}{\mathbf{I}\xspace}

\newcommand{\bL}{\mathbf{L}\xspace}

\newcommand{\bO}{\mathbf{0}\xspace}

\newcommand{\bu}{\mathbf{u}\xspace}
\newcommand{\bv}{\mathbf{v}\xspace}

\usepackage{todonotes}

\title{Review of Cyclic Reduction for Parallel Solution of Hermitian Positive Definite Block-Tridiagonal Linear Systems}

\author{Martin Neuenhofen}

\begin{document}

\maketitle

\begin{abstract}
	Cyclic reduction is a method for the solution of (block-)tridiagonal linear systems. In this note we review the method tailored to hermitian positive definite banded linear systems.
	
	The reviewed method has the following advantages: It is numerically stable without pivoting. It is suitable for parallel computations. In the presented form, it uses fewer computations by exploiting symmetry. Like Cholesky, the reviewed method breaks down when the matrix is not positive definite, offering a robust way for determining positive definiteness.
\end{abstract}

%\tableofcontents

\subsubsection*{Brief Summary} Equations \eqref{eqn:Formulas} give formulas to separate the block-tridiagonal linear system \eqref{eqn:Axb} into two separated block-tridiagonal systems of half the dimension. This can be used for solving the system in a parallel divide-and-conquer approach. The resulting method is described in Algorithm~\ref{algo:proposed}. 

\section{Hermitian Block-Tridiagonal Linear Systems}

We skip the introductory section and literature review on cyclic reduction methods and move directly to the problem.

We consider the linear system
\begin{align}
	\underbrace{\begin{bmatrix}
		\bA_1 	& \bB_1\h 	& 			& 			& 				&  \\
		\bB_1 	& \bA_2 	& \bB_2\h 	& 			& 				&  \\
				& \bB_2 	& \bA_3 	& \bB_3\h 	& 				&  \\
				& 			& \bB_3 	& \ddots 	& \ddots 		&  \\
				& 			&  			& \ddots 	& \ddots 		& \bB_{N-1}\h \\
				& 			& 			& 			& \bB_{N-1} 	& \bA_N
	\end{bmatrix}}_{=:\underline{\bA}} \cdot
	\underbrace{\begin{pmatrix}
		\bx_1\\
		\bx_2\\
		\bx_3\\
		\bx_4\\
		\vdots\\
		\bx_N
	\end{pmatrix}}_{=:\underline{\bx}} = 
	\underbrace{\begin{pmatrix}
	\by_1\\
	\by_2\\
	\by_3\\
	\by_4\\
	\vdots\\
	\by_N
	\end{pmatrix}}_{=:\underline{\by}}\,, \label{eqn:Axb}
\end{align}
where the matrices
\begin{subequations}
	\label{eqn:Data}
	\begin{align}
		\bA_j &\in \C^{m \times m} \quad & \forall\,j&=1,...,N\\
		\bB_j &\in \C^{m \times m} \quad & \forall\,j&=1,...,(N-1)\\
		\by_j &\in \C^{m \times k} \quad & \forall\,j&=1,...,N
	\end{align}
\end{subequations}
are given. $m,k \in \N$ are integers that shall be much smaller than $N \in \N$. The matrix
$$ 	\underline{\bA}\in \C^{(N \cdot m) \times (N \cdot m)} 	$$
shall be positive definite. Sought are the matrices
\begin{align}
	\bx_j &\in \C^{m \times k} \quad & \forall\,j&=1,...,N\,.
\end{align}
We assume $N \in 2 \cdot \N$.

\section{Cyclic Reduction Method for Hermitian Positive Definite Systems}

\paragraph{Derivation of the Method}
The idea of cyclic reduction is based on a special reordering. Putting the elements of $\underline{\bA}$ in the order
\begin{align}
	1,3,5,7,\dots,(N-1)\,,\,2,4,6,8,\dots,N\,,
\end{align}
we obtain the rearranged matrix
\begin{align*}
	\left[
	\begin{array}{ccccc|ccccc}
		\bA_1 	& 			& 			& 			& 			& \bB_1\h 	& 			& 			& 			& 				\\
				& \bA_3 	& 			& 			& 			& \bB_2 	& \bB_3\h 	& 			& 			& 				\\
				& 			& \bA_5 	& 			& 			& 			& \bB_4 	& \bB_5\h	& 			& 				\\
				& 			& 			& \ddots	& 			& 			& 			& \ddots 	& \ddots 	& 	 			\\
				& 			& 			& 			& \bA_{N-1}&			& 			& 			& \bB_{N-2} & \bB_{N-1}\h 	\\
		\hline
		\bB_1 	& \bB_2\h	& 			& 			& 			& \bA_2 	& 			& 			& 			& 				\\
		 		& \bB_3		& \bB_4\h	& 			& 			&  			& \bA_4		& 			& 			& 				\\
		 		& 			& \bB_5		& \ddots	& 			&  			& 			& \bA_6		& 			& 				\\
		 		& 			& 			& \ddots	& \bB_{N-2}\h&  		& 			& 			& \ddots	& 				\\
		 		& 			& 			& 			& \bB_{N-1}&  			& 			& 			& 			& \bA_N				
	\end{array} \right]\,.
\end{align*}
The four sub-blocks in the above matrix we denote with $\underline{\bD}_1,\underline{\bD}_2,\underline{\bC} \in \C^{(N/2\cdot m)\times(N/2 \cdot m)}$\,. I.e., for the above matrix we use the compact writing
\begin{align*}
	\left[\begin{array}{c|c}
	\underline{\bD}_1 	& \underline{\bC}\h\\
	\hline
	\underline{\bC} 	& \underline{\bD}_2
	\end{array}\right]\,.
\end{align*}
These three matrices $\underline{\bD}_1,\underline{\bD}_2,\underline{\bC}$ can be used to construct two separate linear systems for the solution components $\bx_j$. As a benefit, the two separated linear systems have only half the dimension of the original system. In particular, one of the two systems is only in the odd indices $j=1,3,5,\dots,(N-1)$, while the other one is only in the even indices for $\bx_j$. We clarify this in the following.

To be able to write the systems in compact form, we use the notation of even and odd vectors
\begin{align*}
	\underline{\bx}_o = \begin{pmatrix}
		\bx_1\\
		\bx_3\\
		\bx_5\\
		\vdots\\
		\bx_{N-1}
	\end{pmatrix}\,,\quad
	\underline{\bx}_e = \begin{pmatrix}
	\bx_2\\
	\bx_4\\
	\bx_5\\
	\vdots\\
	\bx_{N}
	\end{pmatrix}\,,\quad
	\underline{\by}_o = \begin{pmatrix}
	\by_1\\
	\by_3\\
	\by_5\\
	\vdots\\
	\by_{N-1}
	\end{pmatrix}\,,\quad
	\underline{\by}_e = \begin{pmatrix}
	\by_2\\
	\by_4\\
	\by_5\\
	\vdots\\
	\by_{N}
	\end{pmatrix}\,.
\end{align*}
The decoupled linear systems for $\underline{\bx}_0$ and $\underline{\bx}_e$ are:
\begin{itemize}
	\item Odd System:
	\begin{align}
		\underbrace{(\underline{\bD}_1 - \underline{\bC}\h \cdot \underline{\bD}_2\inv \cdot \underline{\bC})}_{=:\underline{\bU}} \cdot \underline{\bx}_o = \underbrace{\underline{\by}_o - \underline{\bC}\h \cdot \underline{\bD}_2\inv \cdot \underline{\by}_e}_{=:\underline{\bu}}\label{eqn:odd}
	\end{align}
	\item Even System:
	\begin{align}
		\underbrace{(\underline{\bD}_2 - \underline{\bC} \cdot \underline{\bD}_1\inv \cdot \underline{\bC}\h)}_{=:\underline{
		\bV}} \cdot \underline{\bx}_e = \underbrace{\underline{\by}_e - \underline{\bC} \cdot \underline{\bD}_1\inv \cdot \underline{\by}_o}_{=:\underline{\bv}}\label{eqn:even}
	\end{align}	
\end{itemize}
The odd system has the system matrix $\underline{\bU}$ and right-hand side $\underline{\bu}$. The even system has $\underline{\bV}$ and $\underline{\bv}$. In the following we give formulas for these matrices and vectors.

%As we see from the formulas, $\underline{\bU}$ and $\underline{\bV}$ are the primal and dual Schur complements, respectively. Hence, when $\underline{\bA}$ is hermitian positive definite, then so are $\underline{\bU},\underline{\bV}$.

\paragraph{Formulas for the Odd and Even System}
We now give explicit formulas for the matrices $\underline{\bU},\underline{\bV} \in \C^{(N/2 \cdot m) \times (N/2 \cdot m)}$ and the right-hand sides $\underline{\bu},\underline{\bv} \in \C^{(N/2 \cdot m) \times k}$, that build the odd and even system.

The matrices and vectors have the block-structure below:
\begin{subequations}
	\label{eqn:SplitSystems}
\begin{align}
	\underline{\bU} &= 
	\begin{bmatrix}
		\bU_1 	& \bE_1\h 	& 			& 			& 				&  \\
		\bE_1 	& \bU_2 	& \bE_2\h 	& 			& 				&  \\
		& \bE_2 	& \bU_3 	& \bE_3\h 	& 				&  \\
		& 			& \bE_3 	& \ddots 	& \ddots 		&  \\
		& 			&  			& \ddots 	& \ddots 		& \bE_{N/2-1}\h \\
		& 			& 			& 			& \bE_{N/2-1} 	& \bU_{N/2}
	\end{bmatrix}\,,&\quad
	\underline{\bu} &= \begin{pmatrix}
		\bu_1\\
		\bu_2\\
		\bu_3\\
		\vdots\\
		\bu_{N/2-1}\\
		\bu_{N/2}
	\end{pmatrix}\,;\\[8pt]
	\underline{\bV} &= 
	\begin{bmatrix}
	\bV_1 	& \bF_1\h 	& 			& 			& 				&  \\
	\bF_1 	& \bV_2 	& \bF_2\h 	& 			& 				&  \\
			& \bF_2 	& \bV_3 	& \bF_3\h 	& 				&  \\
			& 			& \bF_3 	& \ddots 	& \ddots 		&  \\
			& 			&  			& \ddots 	& \ddots 		& \bF_{N/2-1}\h \\
			& 			& 			& 			& \bF_{N/2-1} 	& \bV_{N/2}
	\end{bmatrix}\,,&\quad
	\underline{\bv} &= \begin{pmatrix}
		\bv_1\\
		\bv_2\\
		\bv_3\\
		\vdots\\
		\bv_{N/2-1}\\
		\bv_{N/2}
	\end{pmatrix}\,.
\end{align}
\end{subequations}
There appear the inner matrices
\begin{subequations}
	\label{eqn:SplitData}
	\begin{align}
		\bU_j,\bV_j &\in \C^{m \times m} & \quad & \forall\,j=1,2,\dots,(N/2)\\
		\bE_j,\bF_j &\in \C^{m \times m} & \quad & \forall\,j=1,2,\dots,(N/2-1)\\
		\bu_j,\bv_j &\in \C^{m \times k} & \quad & \forall\,j=1,2,\dots,(N/2)\,.
	\end{align}
\end{subequations}
The formulas for the these inner matrices can be found by insertion. We present these formulas below:
\begin{subequations}
	\label{eqn:Formulas}
	\begin{align}
		\bU_j &= \bA_{2 \cdot j - 1} - \bB_{2 \cdot j-2} \cdot \bA_{2 \cdot j -2}\inv \cdot \bB_{2 \cdot j-2}\h - \bB_{2 \cdot j-1}\h \cdot \bA_{2 \cdot j}\inv \cdot \bB_{2 \cdot j - 1}\\
		\bE_j &= -\bB_{2 \cdot j} \cdot \bA_{2 \cdot j}\inv \cdot \bB_{2 \cdot j-1}\\
		\bu_j &= \by_{2 \cdot j-1} - \bB_{2 \cdot j -2} \cdot \bA_{2 \cdot j-2}\inv \cdot \by_{2 \cdot j-2} - \bB_{2 \cdot j - 1}\h \cdot \bA_{2 \cdot j}\inv \cdot \by_{2 \cdot j}\\[8pt]
		\bV_j &= \bA_{2 \cdot j} - \bB_{2 \cdot j-1} \cdot \bA_{2 \cdot j -1}\inv \cdot \bB_{2 \cdot j-1}\h - \bB_{2 \cdot j}\h \cdot \bA_{2 \cdot j+1}\inv \cdot \bB_{2 \cdot j}\\
		\bF_j &= -\bB_{2 \cdot j+1} \cdot \bA_{2 \cdot j+1}\inv \cdot \bB_{2 \cdot j}\\
		\bv_j &= \by_{2 \cdot j} - \bB_{2 \cdot j -1} \cdot \bA_{2 \cdot j-1}\inv \cdot \by_{2 \cdot j-1} - \bB_{2 \cdot j}\h \cdot \bA_{2 \cdot j+1}\inv \cdot \by_{2 \cdot j+1}	
	\end{align}
\end{subequations}
The matrices on the left can be computed in parallel for $j=1,2,3,\dots,(N/2)$\,.

It is important that in the above we have -- for the sake of compact formulas -- used identical equations for all $j$. However, for $j=1$ for instance, the formula for $\bU_1$ accesses the matrix $\bB_0$, that does not exist. To repair this detail, we define the following auxiliary matrices, that appear in the above formulas:
\begin{align*}
	\bB_0 &:= \bO \in \C^{m \times m}\,,&\quad
	\bB_N &:= \bO \in \C^{m \times m}\\
	\bA_0 &:= \bI \in \C^{m \times m}\,,&\quad
	\bA_{N+1} &:= \bI \in \C^{m \times m}
\end{align*}

\paragraph{Formulation of the Method}
With the above ingredients, we are able to describe how cyclic reduction is used to solve the linear system \eqref{eqn:Axb}. The method consists of three phases:

In the first phase, given the data \eqref{eqn:Data}, the method uses the formulas in \eqref{eqn:Formulas} to compute the block matrices \eqref{eqn:SplitData} for the odd system \eqref{eqn:odd} and even system \eqref{eqn:even}.

In the second phase, cyclic reduction is used recursively for solving the odd and even system for their respective solutions $\underline{\bx}_o,\underline{\bx}_e$. I.e., we consider \eqref{eqn:odd} and \eqref{eqn:even} as instances of \eqref{eqn:Axb}. This works excellent for the following reason: We assumed that $\underline{\bA}$ is positive definite. If the assumption holds, then further it holds that $\underline{\bU},\underline{\bV}$ are positive definite as well, because they are Schur complements of $\underline{\bA}$.

In the third phase, we rearrange the vectors $\underline{\bx}_o,\underline{\bx}_e$ into the global solution vector $\underline{\bx}$ of \eqref{eqn:Axb}. This completes the method.

\paragraph{Determination of Positive Definiteness}
If the matrix $\underline{\bA}$ is indeed positive definite, then all recursive calls of cyclic reduction will succeed, until at the bottom there arise positive definite linear systems of size $m \times m$. These shall be solved with Cholesky's method.

If, in contrast to the assumption, the matrix $\underline{\bA}$ is not positive definite, then eventually at least either of the matrices $\underline{\bU},\underline{\bV}$ will not be well-defined, and thus the method will crash. This crash will appear in the form that Cholesky's method determines indefiniteness for one of the matrices $\bA_j$ in \eqref{eqn:Formulas}.

\section{Description of a Parallel Method}

\paragraph{Computing Model}
For our description of the parallel method, we consider a parallel computing system that consists of $N$ computing nodes; each has an index, counted as $1,2,3,...,N$. The nodes are connected in a communication network, that consists of cables. The length each cable equals the distance of the respective nodes in space. Information travels along these cables at a limited speed.

\paragraph{Idea of the Parallel Method}
The idea for a parallel cyclic reduction consists of two ingredients:
\begin{enumerate}
	\item \textbf{Parallel instantiation of odd and even systems:} The equations \eqref{eqn:Formulas} are evaluated in parallel over $j=1,\dots,(N/2)$. If $N$ parallel computing units are given, then the nodes $1,2,\dots,N/2$ can compute the matrices for the odd system, and the nodes with indices $N/2+1,\dots,N$ can compute the matrices for the even system.
	\item \textbf{Parallel recursive solution of the split systems:} The linear systems \eqref{eqn:odd} and \eqref{eqn:even} are solved recursively for $\underline{\bx}_o$ and $\underline{\bx}_e$. To this end, the same method as for \eqref{eqn:Axb} can be used.
\end{enumerate}

A pseudo-code for this procedure is given in Algorithm~\ref{algo:proposed}. The algorithm assumes $N \in 2^\N$\,.

\begin{algorithm}
	\caption{Cyclic Reduction Method}
	\label{algo:proposed}
	\begin{algorithmic}[1]
		\Procedure{Solver}{$\underline{\bA},\underline{\by}$}
		\If{$\dim(\underline{\bA})$ is sufficiently small}
			\State Solve $\underline{\bA} \cdot \underline{\bx} = \underline{\by}$ for $\underline{\bx}$, using Cholesky's method.
			\State \Return $\underline{\bx}$
		\EndIf
		\State \textit{// parallel for-loop for nodes $j=1,2,\dots,(N/2)$}
		\For{$j=1,2,\dots,N/2$}
		\State Compute $\bU_j,\bE_j,\bu_j$ from formulas \eqref{eqn:Formulas}.
		\EndFor
		\State \textit{// concurrent parallel for-loop for nodes $j=(N/2)+1,(N/2)+2,\dots,N$}
		\For{$j=1,2,\dots,N/2$}
		\State Compute $\bV_j,\bF_j,\bv_j$ from formulas \eqref{eqn:Formulas}.
		\EndFor
		\State \textit{// The even and odd systems $\underline{\bU},\underline{\bV}$ are given as $\bU_j,\bE_j,\bV_j,\bF_j,\bu_j,\bv_j$ according to \eqref{eqn:SplitSystems}}
		\State $\underline{\bx}_o := $\textsc{Solver}($\underline{\bU},\underline{\bu}$)\quad\textit{// solve on nodes $j=1,2,\dots,(N/2)$}
		\State $\underline{\bx}_e := $\textsc{Solver}($\underline{\bV},\underline{\bv}$)\quad\textit{// solve on nodes $j=(N/2)+1,(N/2)+2,\dots,N$}
		\State Reassemble $\underline{\bx}$ from $\underline{\bx}_o,\underline{\bx}_e$\,.
		\State \Return $\underline{\bx}$
		\EndProcedure
	\end{algorithmic}
\end{algorithm}

\paragraph{Time-Complexity and Communication}
The pure time-complexity of parallel computations is in
$$	\cO(\,\log(N)\cdot m^3\,)\,. 	$$
But, there is an additional cost for communication of data between the computing nodes. This cost is discussed below.

\largeparbreak
The method has significant non-local communication. We describe this in the following, where $j$ ranges from $1$ to $N/2$:

For the computation of the odd system, the node $j$ receives data from nodes $2 \cdot j-2$, $2 \cdot j-1$, and $2 \cdot j$\,. Analogously, for the computation of the even system, the node $N/2+j$ receives data from nodes $2 \cdot j-1$, $2 \cdot j$, and $2 \cdot j+1$\,.

Depending on the communication network of the parallel computing system, this non-local communication can be very expensive. In fact, it may result in a time complexity that exceeds the time spent in actual parallel computations.

Non-local communications are not unavoidable in parallel methods for solving block-tridiagonal linear systems.
In \cite{myblocktridiagLGS} we provide an accurate description of a method (with directly implementable code for an MPI cluster) that works only via local communications. The communication network used for this method is a binary tree -- the cheapest possible way for connecting all nodes.

\section{Numerical stability}
If $\underline{\bA}$ is positive definite, then: The matrices $\underline{\bD}_1,\underline{\bD}_2$ are regular, and the matrices $\underline{\bQ}_o$ and $\underline{\bQ}_e$ are positive definite.
If $\underline{\bA}$ is positive definite, and Cholesky decomposition with forward and backward substitution is used to compute formulas \eqref{eqn:Formulas}, then: The matrices $\underline{\bQ}_o$ and $\underline{\bQ}_e$ are bit-wise hermitian, and the algorithm is numerically stable. In particular, the substitution must be applied as follows:

Given a formula like
\begin{align*}
	\bZ := \bG\h \cdot \bM\inv \cdot \bG
\end{align*}
we must use the Cholesky decomposition $M = \bL \cdot \bL\h$ in the following way:
\begin{algorithmic}[1]
	\State $\bL :=$ \textsc{Chol}($\bM$)
	\State $\tilde{\bG} := \bL \backslash \bG$
	\State \Return $\bZ := \tilde{\bG}\h \cdot \tilde{\bG}\h$
\end{algorithmic}

Using the cyclic reduction in a recursive way, or solving the split systems with an arbitrary forward stable numerical algorithm, results in an overall method that is numerically stable.
\largeparbreak

The described algorithm can also be used to determine whether a matrix $\underline{\bA}$ at hand is positive definite. Namely, this is the case if and only if all split subsystems are positive definite. Hence, we can stop the recursive splitting when the submatrices are small enough, and then check, e.g. with Cholesky, for the positive definiteness of $\underline{\bU},\underline{\bV}$.

\largeparbreak

\FloatBarrier

%\nocite{*}
\bibliography{CyclicReduc_bib}
\bibliographystyle{plain}

\end{document}